\newcounter{firstpage} 
\newcounter{finalpage} 
\documentclass[11pt,twoside]{article}
\usepackage[tbtags]{amsmath}
\usepackage{amsthm,amsfonts,amssymb,verbatim,cite}
\makeatletter
\def\theorem@style{\thm@style}
\def\theorembodyfont{\thm@bodyfont}
\def\theoremheadfont{\thm@headfont}
\def\theoremnotefont{\thm@notefont}
\def\theoremheadpunct{\thm@headpunct}
\def\theorempreskipamount{\thm@preskip}
\def\theorempostskipamount{\thm@postskip}
\def\theoremindent{\thm@indent}
\makeatother
\makeatletter
\let\over\@@over
\let\atop\@@atop
\let\@savedendocument\enddocument
\def\enddocument{%
\clearpage
\setcounter{finalpage}{\thepage} \addtocounter{finalpage}{-1}
\immediate\write\@auxout{\string\setcounter{finalpage}{\arabic{finalpage}}}
\markboth{}{}\cleardoublepage
\@savedendocument}
\makeatother
\advance\baselineskip 0pt plus 1.5pt
\advance\parskip 0pt plus 2pt
\makeatletter\@addtoreset{equation}{section}\makeatother

\newcommand{\TitAuthAffilArch}[4]{%
\leftline{\copyright~ 1998 International Press}%
\leftline{Adv. Theor. Math. Phys. {\bf 2} (1998) \thepage--\thefinalpage}%
\vspace{2cm}
\renewcommand{\thefootnote}{\relax}
\begin{center}
\LARGE\bf #1%
\footnote{e-print archive: \texttt{http://xxx.lanl.gov/abs/#4}}
\end{center}
\renewcommand{\thefootnote}{\arabic{footnote}}\setcounter{footnote}{0}
\begin{center}
{\bf #2}\\
\vspace{0.5cm}
#3
\end{center}
\vspace{1cm}
}

\theoremstyle{plain}
\newtheorem{thm}{Theorem}[section]

\theoremstyle{definition}

\theoremstyle{remark}

\newcommand{\Tr}{\mathop{\mathgroup\symoperators Tr}\nolimits}

\makeatletter
\def\@begintheorem#1#2[#3]{%
\item[\normalfont %
\hskip\labelsep
\the\theoremheadfont
\theoremindent
\@ifempty{#1}{\let\thmname\@gobble}{\let\thmname\@iden}%
\@ifempty{#2}{\let\thmnumber\@gobble}{\let\thmnumber\@iden}%
\let\thmnote\@gobble
\thm@swap\swappedhead\thmhead{#1}{#2}{#3}%
]%
\hskip-\labelsep
{\normalfont \the\theoremheadfont
\@ifempty{#3}{\let\thmnote\@gobble}{\let\thmnote\@iden}%
\thmnote{\hbox{ }{\the\theoremnotefont(#3)}}\the\theoremheadpunct
\hskip\labelsep}%
\thmheadnl %
\ignorespaces}
\makeatother

\newcommand{\eqnref}[1]{(\ref{#1})}

\makeatletter
\renewcommand\section{\@startsection {section}{1}{\z@}%
{-3.5ex \@plus -1ex \@minus -.2ex}%
{2.3ex \@plus.2ex}%
{\normalfont\Large\bfseries\boldmath}}
\renewcommand\subsection{\@startsection{subsection}{2}{\z@}%
{-3.25ex\@plus -1ex \@minus -.2ex}%
{1.5ex \@plus .2ex}%
{\normalfont\large\bfseries\boldmath}}
\renewcommand\subsubsection{\@startsection{subsubsection}{3}{\z@}%
{-3.25ex\@plus -1ex \@minus -.2ex}%
{1.5ex \@plus .2ex}%
{\normalfont\normalsize\bfseries\boldmath}}
\makeatother
\makeatother \usepackage{epsf,graphicx}
\makeatletter 
\newcommand{\Dsl}{\,\raise.15ex\hbox{/}\mkern-13.5mu D} 

\providecommand{\eqref}[1]{(\ref{#1})}
\providecommand{\@gobble}[1]{}
\newcommand{\HMlref}[2]{\unskip\newcommand{#1}{\cite{\expandafter\@gobble\string#1}}}

\DeclareRobustCommand\HMrefs{%
  \@ifnextchar [{\@tempswatrue\HM@citex}{\@tempswafalse\HM@citex[]}}
\def\HM@citex[#1]#2{%
  \let\@citea\@empty
  \@cite{\@for\@citeb:=#2\do
    {\@citea\def\@citea{,\penalty\@m\ }%
      {\def\cite#1{#1}\xdef\@citeb{\@citeb}}%
      \edef\@citeb{\expandafter\@firstofone\@citeb}%
     \if@filesw\immediate\write\@auxout{\string\citation{\@citeb}}\fi
     \@ifundefined{b@\@citeb}{\mbox{\reset@font\bfseries ?}%
       \G@refundefinedtrue
       \@latex@warning
         {Citation `\@citeb' on page \thepage \space undefined}}%
       {\hbox{\csname b@\@citeb\endcsname}}}}{#1}}

\newcommand{\enlabel}[1]{%
   \expandafter\gdef\csname #1\endcsname{\eqref{#1}}%
   \label{#1}
}

\newcommand{\HMeqn}[2]{%
   \gdef#1{\eqref{\expandafter\@gobble\string#1}}
   \begin{equation}\relax
   \expandafter\label\expandafter{\expandafter\@gobble\string#1}
   #2 \end{equation}
}

\makeatother 
\advance\baselineskip 0pt plus 1pt
\advance\parskip 0pt plus 2pt

\DeclareMathOperator{\SU}{SU} 
 
\DeclareMathOperator{\sll}{sl} 

\DeclareMathOperator{\End}{End}
\newcommand\Xcite\cite

\def\epsfsize#1#2{0.4#1} 

\newcommand{\Uqsl}{{\text{U}_{q}\sll}} 
 
\newcommand{\C}{\mathbb{C}}

\newcommand{\R}{\mathbb{R}}

\renewcommand{\d}{{\text{d}}}

\begin{document}

\TitAuthAffilArch{The Classical Evaluation of Relativistic Spin Networks}%
{John W. Barrett\footnote{e-mail: \texttt{jwb@maths.nott.ac.uk}}}%
{Department of Mathematics \\
 University of Nottingham \\
 University Park \\
 Nottingham, NG7 2RD, UK}{math.QA/9803063}

\begin{abstract}
The evaluation of a relativistic spin network for the classical case of the
Lie group is given by an integral formula over copies of
$\SU(2)$. For the graph determined by a 4-simplex this gives the evaluation
as an integral over a space of geometries for a 4-simplex.
\end{abstract}

\section{Introduction}

 A relativistic spin network is a graph embedded in $\R^3$ with a non-negative integer, the spin, labelling each edge of the graph. An evaluation of relativistic spin networks was defined in \cite{Barr:Cr}. The evaluation gives an invariant of the isotopy class of the labelled embedded graph. The invariant is a rational function, with integer coefficients, of a parameter $A$, which is the parameter for the deformation of $\sll(2)$.

In \cite{Barr:Cr}, we discussed the case of 4-valent graphs, giving a canonical formula. This was generalised to arbitrary graphs in \cite{Yett}, which gave the complete normalisations and established some important properties of the invariant.

In this paper, an integral formula is presented for the evaluation in the classical case, which is when the parameter $A$ is specialised to $\pm1$. In this case, the invariant is a rational number which depends only on the graph and not on the embedding.
The first section presents the formula for this invariant as an integral over several copies of the Lie group $\SU(2)$. 

The second section relates this construction to the original combinatorial definition of the relativistic spin network evaluation, giving a proof that it is the specialisation of the combinatorial definition to $A=\pm1$. Some identities which are essentially some simple examples of this theorem are given in \cite{Tan:Tan}. 

The third section gives some remarks on the geometrical interpretation of the new integral formula in terms of four-dimensional Euclidean geometry. In particular, for the graph determined by a 4-simplex the evaluation is an integral over a space of geometries for a 4-simplex.

The evaluation of \cite{Barr:Cr} is defined using the representation theory of the quantum group $\Uqsl(2)$, with $q=A^2$. The classical value $q=1$ corresponds to the classical complex Lie algebra $\sll(2)$, and the spins label the irreducible  representations. Spin $n$ corresponds to the representation of dimension $n+1$.  These can also be regarded as the irreducible representations of the compact Lie group $\SU(2)$, and so the classical invariant can be written in terms of this group.

\section{The Graph Invariant}

The notation for the irreducible representation of $\SU(2)$ of spin $n$ is as follows. The vector space is denoted $V$ and the representation is $\rho\colon\SU(2)\to\End(V)$, where $\End(V)$ denotes the vector space of $\C$-linear maps $V\to V$.  

A graph is a 1-dimensional cell complex. This is a finite set of points called vertices, together with a finite set of 1-cells, the edges, with the ends of the edges identified with vertices. It is possible for an several edges to join a pair of vertices (multiple edges) or for a edge to start and end at the same vertex (a loop). A vertex is called $n$-valent if the number of edges meeting it (its degree) is $n$.

A variable $h_k\in\SU(2)$ is assigned to each vertex $k$. Then a weight in $\R$ is defined for each edge $e$ of the graph, with its vertices denoted by $e(0)$ and $e(1)$, as
\begin{equation}
(-1)^{n_e} \Tr \rho_{e}\left(h_{e(0)}^{\strut} h_{e(1)}^{-1}\right) 
\label{tag1}
\end{equation}
In this formula, $n_e$ is the spin of edge $e$ and $\rho_e$ the corresponding representation. This weight does not depend on the order of the vertices of the edge.
 
The invariant $I\in \R$ is defined by taking the product of all the weights over the set of edges $E$, and integrating over each copy of $\SU(2)$ for each vertex. The integration is done using the Haar measure on $\SU(2)$, normalised to total volume 1.  The formula is
$$ I=(-1)^{\sum_E n_{e}} \int_{h\in SU(2)\times\ldots\times\SU(2)} \prod_{e\in E} \Tr \rho_{e}\left(h_{e(0)}^{\strut} h_{e(1)}^{-1}\right).$$
The factors of $-1$ are not an essential feature of this invariant, but have been included to agree with the relativistic spin network definition for general $A$. 

This formula clearly generalises to other compact Lie groups which have the property that every element is conjugate to its inverse.

The integration over the group ensures that the invariant is zero unless the tensor product of the representations at each vertex contains the trivial representation. This implies certain restrictions on the spins on the edges incident to one vertex. One restriction is that the sum of the spins is an even integer. The other restriction is that it is possible to assign a vector in $\R^3$ to each edge meeting a given vertex, with lengths given by the corresponding spins, such that the vectors sum to zero. For example, with a trivalent vertex this gives the familiar triangle inequalities which occur in the Clebsch-Gordan series for $\SU(2)$. For a bivalent vertex it implies that the spins on the two edges must be equal, and for a monovalent vertex, the spin must be 0. In this case, the edge can be removed without altering the invariant.

An edge which starts and ends at the same vertex contributes a factor of $$(-1)^{n}\Tr\rho(1)=(-1)^{n}(n+1),$$
known as the quantum dimension of the representation $\rho$. This loop can be removed from the graph if the invariant is multiplied by this factor.

Let $\Gamma'$ be a graph obtained from $\Gamma$ by the contraction of an edge joining distinct vertices. For example, a part of the graph  
$${\def\epsfsize#1#2{0.2#1}  
\hspace{20mm}
\vcenter{
\epsfbox{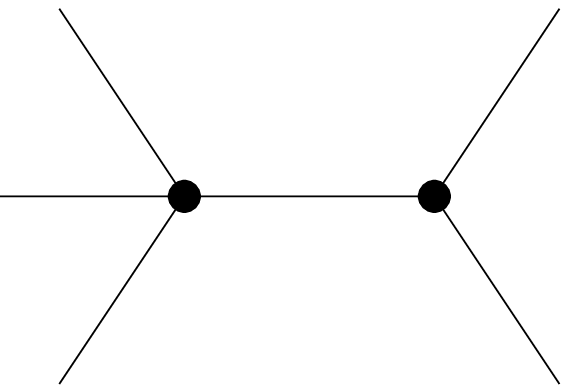}}
\quad\text{is replaced by }\quad
\vcenter{
\epsfbox{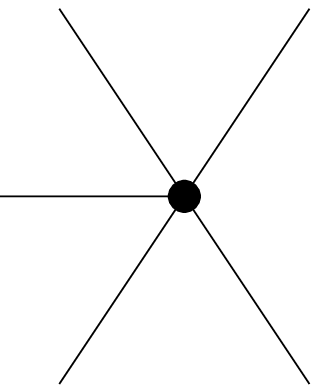}}
}$$
in a contraction. The invariant for $\Gamma'$ with the edges labelled by spins is obtained from the invariant for $\Gamma$ with the corresponding edges labelled the same, by summing over the spin $n$ on the contracted edge,
$$ I_{\Gamma'}=\sum_{n=0}^{\infty}
 (-1)^{n}(n+1) I_\Gamma.$$
Only a finite number of terms in this sum are not zero, as there is an upper bound on the spin $n$ meeting the restrictions given above.

The proof of this contraction formula follows from the character formula
$$\sum_{n=0}^{\infty} (n+1)\int_{(g,h)\in\SU(2)\times\SU(2)} f(g,h)\Tr\rho\left(g h^{-1}\right) =
\int_{h\in\SU(2)}f(h,h)$$
for any continuous function $f$.

A corollary of the contraction formula is that any vertex can be expanded as a tree in any fashion, with a sum over the spins on the edges of a tree, as noted in the more general $q$ context in \cite{Yett}.

The contraction formula does not apply to a loop. Summing over the spins on a loop would give
$$ I_{\Gamma'}\sum_{n=0}^{\infty}(n+1)^2$$
which does not converge.

Multiple edges can be replaced in the formula by a sum over spins on a single edge connecting the two vertices. This follows since
$$\prod_e \Tr \rho_e(g)=\Tr\rho_{\otimes e}(g)$$
and the tensor product can be replaced by the sum of irreducible representations, giving a sum of spins on a single edge.

\section{Relativistic Spin Network Evaluations}

The relativistic spin network evaluation is defined in terms of the ordinary $SU(2)$ spin network evaluation, as defined by \cite{Pen} for the case $A=-1$, and generalised to the $q$-deformed case, with arbitrary $A$ in \cite{Kauff}. The specific normalisation used here is given in \cite{Kauff:Lins}. The term spin network (without the adjective relativistic) will refer to the original spin networks considered by Penrose and Kauffman.

The spin network is a trivalent ribbon graph in $\R^3$ with spin labels on the edges. Ribbon means that the graph is embedded in an oriented surface with boundary, which can be regarded as a thickening of the original graph \cite{Resh:Tur}. This allows the spin network evaluation to take into account the cyclic order of edges at a vertex and the framing of the edges. This data can alternatively be specified by a plane projection of the graph, as was originally envisaged by Penrose. 

Denote the spin network evaluation by $N(A)$. The definition of the relativistic spin network evaluation for a trivalent graph is
$$ R(A)={ N(A) N(A^{-1})\over \prod_{\text{vertices}} \theta},$$
where the normalisation factor of $\theta$ at each vertex is the evaluation of the $\theta$ graph formed by closing two copies of the spin network vertex to a closed graph.
Alternatively, one can work with the spin network vertices normalised so that each $\theta=1$.

The definition given in \cite{Barr:Cr} is the product of the two spin network evaluations, where the second one is obtained from the first by interchanging over- and under-crossings in a planar projection. However this is the same as replacing $A$ with $A^{-1}$. Now $R(A)=R(-A)$ due to the fact that $N(A)$ changes by a power of $-1$ for each crossing point, and these cancel in $R(A)$. Thus $R$ depends only on $q=A^2$.

The definition of the spin network evaluation for arbitrary graphs is given in \cite{Yett}. The definition is given by applying the contraction formula to a trivalent graph, as in the previous section. In this more general setting, the weight for the spin on the contracted edge is the quantum dimension of the corresponding representation of $\Uqsl(2)$.\footnote{In part of the literature, $q$ is taken to be $A^4$, which can cause confusion. The conventions are explained in \cite{Saw}.} Yetter shows that this definition does not depend on the trivalent tree which is chosen to expand a vertex (of higher degree). When $q$ is a non-trivial root of unity, the range of summation in the contraction formula has to be restricted to a finite range. Namely, if $r$ is the smallest integer such that $q^{2r}=1$, and $r>1$, then the sums are over spins $0\le n\le r-2$. The required formulae are developed for the case $q=e^{i\pi/r}$ in \cite{Kauff:Lins}. 

The spin network evaluation is defined for planar projections of ribbon graphs with free ends. A graph with free ends is more general than the graphs defined previously, in that the free ends of edges are ends which are not identified with vertices.  A plane projection of an embedded graph with free ends is required to lie in $\R\times [0,1]$ with the free ends on the boundary. These 
form a category with objects given by the free ends, and can be composed in an obvious way by connecting two edges with free ends that meet to form one edge. Then the evaluation, for a specific $A\in\C$, is a monoidal functor to a category $L_A$ in which the set of endomorphisms of the empty object (no free ends) is the coefficient field $\C$. This construction was given in \cite{Resh:Tur} for the case when the target category is the category of vector spaces, which works when $q$ is not a non-trivial root of unity $(r>1)$. For the non-trivial root of unity cases, the construction is a quotient 
\cite{Barr:WestI,Barr:WestII}.

The relativistic spin network evaluation with free ends can be described in terms of the spin network evaluation with free ends. In the monoidal category $C$ of trivalent tangle diagrams let $m^\dag$ be the
morphism $m$ reflected in a horizontal line. If the spin network
evaluation of $m$ in the linear category $L_A$ is a map $l(m)\colon V\to W$, then the relativistic spin network evaluation of $m$ is the linear map $\End(V)\to\End(W)$ given by $\alpha\mapsto l(m^\dag)\circ \alpha \circ l(m)$. This evaluation is a monoidal functor from $C$ to the category of vector spaces.
  
In the following equations, the relativistic spin network appears on the left and spin networks on the right. The $1/\theta$ normalisation factors have been omitted.  The equations show simple cases of the evaluation where $V$ is taken to be the trivial object, and $\alpha=1_V$.

The definition of the relativistic vertex is
\begin{equation}
\vcenter{\epsfbox{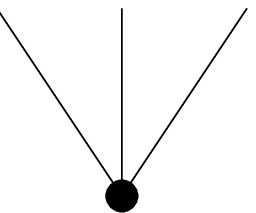}}=\vcenter{\epsfbox{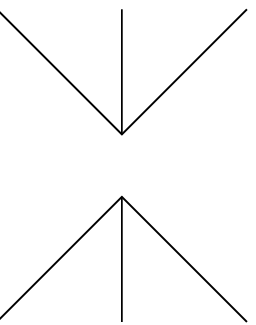}}
\label{tag2}
\end{equation}
In the right hand term the top vertex is labelled the same as the left hand relativistic vertex, and the bottom vertex is the top one reflected in a horizontal line.

The braiding is
$$\vcenter{\epsfbox{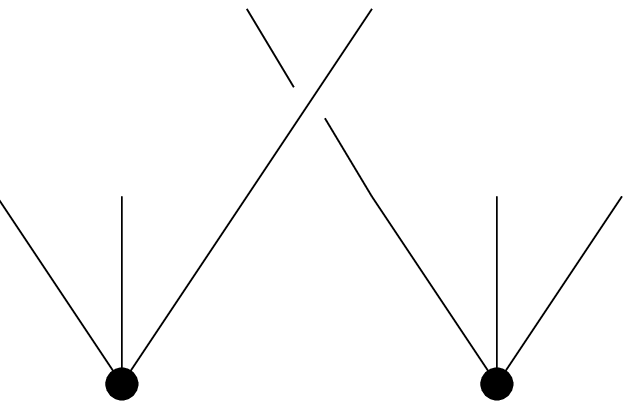}}=\vcenter{\epsfbox{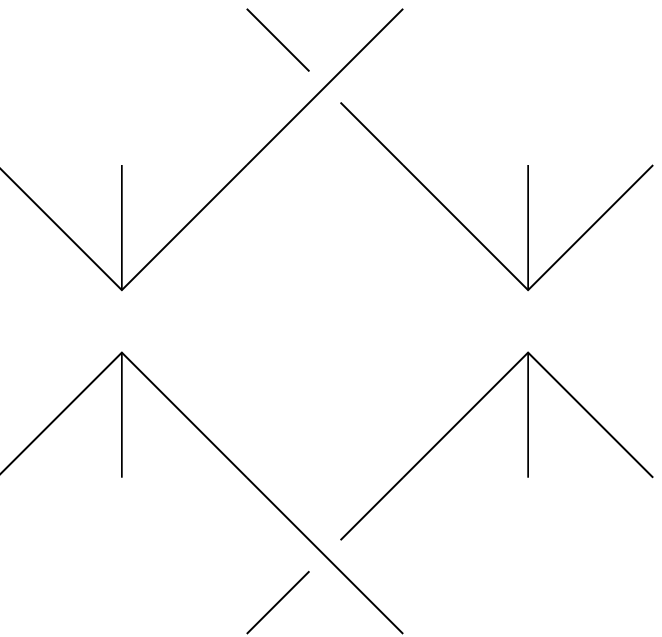}}$$
and the inner product, which joins free ends with the same spin,
\begin{equation}\
\vcenter{\epsfbox{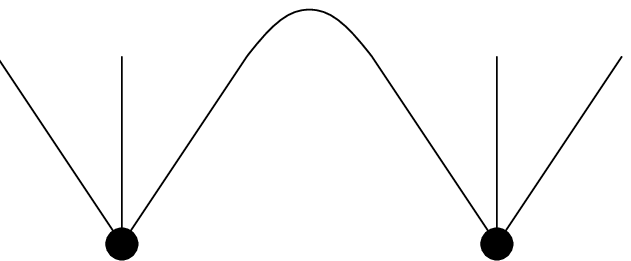}}=\vcenter{\epsfbox{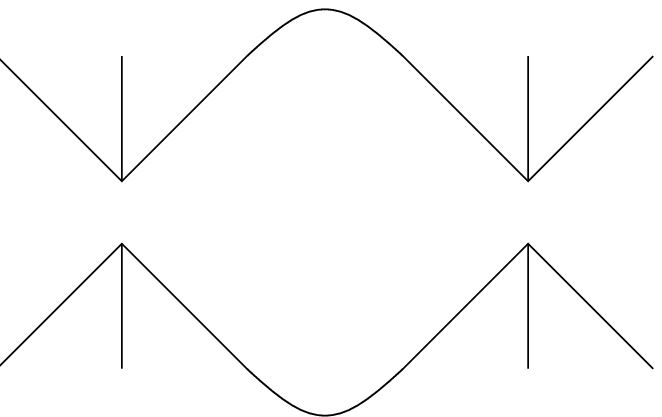}}\label{tag3}
\end{equation}

These cases are sufficient for the construction of any graph from these elementary vertices, braidings, maxima and minima. This is because any graph can be presented with the vertices and minima at the bottom of the diagram, as in \eqnref{tag2}. 

Now the main result can be presented.

\begin{thm} 
The relativistic spin network evaluation specialises at $A=\pm1$ to give the 
graph invariant $I$ defined in the previous section. 
\label{Theorem:3.1}
\end{thm}

\begin{proof}
For $A=\pm1$, the spin networks are evaluated in the category of representations of $\SU(2)$, with morphisms $\SU(2)$-invariant linear maps \cite{Mouss,Barr}. This induces the following representation for the relativistic spin networks. A free end labelled $n$ is sent to the vector space $\End(V)$, where $V$ is the corresponding vector space for the representation. The trivalent vertex \eqnref{tag2}
is the projector in $\End(V_1)\otimes\End(V_2)\otimes\End(V_3)$ onto the subspace of linear maps which commute with the action of $\SU(2)$. This projector is equal to
\begin{equation}
\int_{\SU(2)}\d g\quad \rho_1(g)\otimes\rho_2(g)\otimes\rho_3(g).
\label{tag4}
\end{equation}
The inner product is the map $\End(V)\otimes\End(V)\to\C$ given by
$$\langle\alpha,\beta\rangle =\epsilon\circ(\alpha\otimes \beta)\circ\eta,$$
where $\epsilon\colon \C\to V\otimes V$ is the `minimum' \;\epsfbox{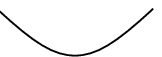}\; of \eqnref{tag3} and $\eta\colon  V\otimes V\to\C$ is the `maximum' \;\epsfbox{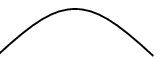}\,. This is exactly the inner product used in \eqnref{tag3}. Now consider the special case where $\alpha$ and $\beta$ are the action of elements $g, h\in\SU(2)$. Since the minima and maxima are $\SU(2)$-invariant, 
 $$\langle\rho(g),\rho(h)\rangle=
\epsilon\circ\left(1\otimes\rho(g)^{-1}\rho(h) \right)\circ\eta=
(-1)^{n}\Tr\rho\left(h^{-1}g\right)$$
 This formula is the same as \eqnref{tag1}. Finally, the braiding is just the twist map $\alpha\otimes\beta\mapsto\beta\otimes\alpha$.
Combining the vertex \eqnref{tag4} with the inner product gives the invariant $I$ for a closed trivalent graph. Since the contraction formula is the same for both invariants, they coincide for graphs with vertices of any valence.
\end{proof}

\section{Geometrical Interpretation}
 Using the identification $\SU(2)\cong S^3$, each variable $h\in\SU(2)$ at a vertex can be regarded as a unit vector in $\R^4$. Then the character \eqnref{tag1} is a function of the angle $\phi$ between the two vectors defined by $h_{e(0)}$ and $h_{e(1)}$,
\begin{equation}
\Tr \rho\left(h_{e(0)}^{\strut} h_{e(1)}^{-1}\right)={\sin(n+1)\phi\over\sin\phi}. 
\label{tag5}
\end{equation}

Consider a compact simplicial 3-manifold. The dual 1-skeleton $\Gamma$ is a graph which has a vertex for every tetrahedron in the manifold, and an edge for every triangle, with free ends for the boundary triangles. The unit vector $h\in S^3$ at each vertex of the graph can be regarded as the normal to a (3-dimensional) hyperplane, now associated to the tetrahedron of the 3-manifold.  In the non-trivial case when the sum of the spins at a vertex of $\Gamma$ is even, the product of the weights \eqnref{tag5} is unchanged if $h$ is replaced by $-h$ at any vertex. Thus the integrand does not register the orientation of the hyperplane.

In \cite{Barr:Cr}, we discussed the example of the boundary of a 4-simplex. In this case the graph $\Gamma$ is the complete graph on five vertices. We gave an argument, based on quantization, that the relativistic spin network evaluation is related to the metric geometry of a 4-simplex embedded in 4-dimensional Euclidean space $\R^4$. In the integral formula presented here, the five hyperplanes corresponding to the five vertices of $\Gamma$ determine a geometric 4-simplex in $\R^4$ uniquely up to translation and an overall scaling (which can be positive or negative).  

In this way, the weights \eqnref{tag5} depend on the ten dihedral angles of the 4-simplex. These ten dihedral angles characterise the metric geometry of a 4-simplex uniquely up to isometry and scaling. The invariant $I$ is thus determined by averaging over these geometries.

\end{document}